\documentclass[12pt,reqno]{amsart}
\usepackage{amsmath, amsfonts, amssymb, amsthm}
\def\Z{\mathbb Z}
\def\C{\mathbb C}

\def\Q{\mathbb Q}
\def\R{\mathbb R}

\def\cP{\mathcal P}

\def\1{{\bf 1}}

\def\qbinom #1#2#3{{\genfrac{[}{]}{0pt}{}{#1}{#2}}_{#3}}
\theoremstyle{plain}
\newtheorem{theorem}{Theorem}[section]

\theoremstyle{definition}

\theoremstyle{remark}
\newtheorem*{Rem}{Remark}

\numberwithin{equation}{section}

\begin{document}

\title{Some symmetric $q$-congruences}

\begin{abstract}
We prove some symmetric $q$-congruences.
\end{abstract}
\author{He-Xia Ni} 
\address{Department of Mathematics, Nanjing University, Nanjing 210093,
People's Republic of China}
\email{nihexia@yeah.net}
\author{Hao Pan}
\address{Department of Mathematics, Nanjing University, Nanjing 210093,
People's Republic of China}
\email{haopan79@zoho.com}

\maketitle

\section{Introduction}
\setcounter{theorem}{0}\setcounter{lemma}{0}\setcounter{equation}{0}

Suppose that $f_0,f_1,f_2,\ldots$ are a sequence of integers. Let
$$
\hat{f}_k(q)=\sum_{j=0}^k(-1)^j\binom{k}jf_j.
$$
In \cite[Theorem 2.4]{Sun14}, Sun proved the following symmetric congruence:
\begin{equation}\label{suncong}
\sum_{k=0}^{p-1}\binom{\alpha}{k}\binom{-1-\alpha}{k}f_k\equiv (-1)^{\langle\alpha\rangle_p}
\sum_{k=0}^{p-1}\binom{\alpha}{k}\binom{-1-\alpha}{k}\hat{f}_k\pmod{p^2},
\end{equation}
where $p$ is an odd prime, $\alpha\in\Q$ is $p$-integral and $\langle\alpha\rangle_p$ is the least non-negative residue of $\alpha$ modulo $p$. With the help of (\ref{suncong}), Sun obtained many interesting congruences modulo $p$.

In this paper, we shall give a $q$-analogue of (\ref{suncong}). Define
$$
(x;q)_n=\begin{cases}
(1-x)(1-xq)\cdots(1-xq^{n-1}),&\text{if }n\geq 1,\\
1,&\text{if }n=0.
\end{cases}
$$
Define the $q$-binomial coefficient
$$
\qbinom{\alpha}{k}q=\frac{(q^{\alpha-k+1};q)_{k}}{(q;q)_k}.
$$
In particular, we set $\qbinom{\alpha}{k}q=0$ if $k<0$.
Suppose that $f_1(q),f_2(q),\ldots$ are a sequence of polynomials in $q$ with integral coefficients. Let
$$
\hat{f}_k(q)=\sum_{j=0}^k(-1)^jq^{\binom{j+1}{2}}\qbinom{k}{j}qf_j(q)
$$
for $0\leq k\leq n$.
\begin{theorem}\label{qsymthm} Let $n\geq 2$ and $d\geq 1$ with $(n,d)=1$. 
Suppose that $f_1(q),\ldots,f_n(q)\in\Z[q]$ and $r$ is an integer.
Then if $n$ is odd,
\begin{align}\label{qsymodd}
&q^{d\binom{a+1}2+\frac{(ad+r)(n-1-2a)}{2}}\sum_{k=0}^{n-1}
\frac{(q^r;q^d)_k(q^{d-r};q^d)_k}{(q^d;q^d)_k^2}\cdot q^{dk}f_k(q^d)\notag\\
\equiv&
(-1)^a\sum_{k=0}^{n-1}\frac{(q^r;q^d)_k(q^{d-r};q^d)_k}{(q^d;q^d)_k^2}\cdot q^{dk}\hat{f}_k(q^d)\pmod{\Phi_n(q)^2},
\end{align}
where $a=\langle-r/d\rangle_n$ and
$\Phi_n$ is the $n$-th cyclotomic polynomial.
And when $n$ is even,
\begin{align}\label{qsymeven}
&q^{d\binom{a+1}2+\frac{(ad+r)(n-1-2a)}{2}}\sum_{k=0}^{n-1}
\frac{(q^r;q^d)_k(q^{d-r};q^d)_k}{(q^d;q^d)_k^2}\cdot q^{dk}f_k(q^d)\notag\\
\equiv&
(-1)^{a+\frac{ad+r}{n}}\sum_{k=0}^{n-1}\frac{(q^r;q^d)_k(q^{d-r};q^d)_k}{(q^d;q^d)_k^2}\cdot q^{dk}\hat{f}_k(q^d)\pmod{\Phi_n(q)^2}.
\end{align}
\end{theorem}
Now for odd $n\geq 3$, replacing $q$ by $q^{-1}$ in (\ref{qsymodd}) and noting that
$$
\Phi_n(q^{-1})=q^{-\phi(n)}\Phi_n(q)
$$
where $\phi$ is the Euler totient function, we get
\begin{align*}
&q^{-d\binom{a+1}2-\frac{(ad+r)(n-1-2a)}{2}}\sum_{k=0}^{n-1}
\frac{(q^{-r};q^{-d})_k(q^{r-d};q^{-d})_k}{(q^{-d};q^{-d})_k^2}\cdot q^{-dk}f_k(q^{-d})\\
\equiv&
(-1)^a\sum_{k=0}^{n-1}\frac{(q^{-r};q^{-d})_k(q^{r-d};q^{-d})_k}{(q^{-d};q^{-d})_k^2}\cdot q^{-dk}\hat{f}_k(q^{-d})\pmod{\Phi_n(q)^2}.
\end{align*}
Notice that
$$
\hat{f}_k(q^{-1})=
\sum_{j=0}^k(-1)^jq^{-\binom{j+1}{2}}\qbinom{k}{j}{q^{-1}}f_j(q^{-1})
=\sum_{j=0}^k(-1)^jq^{\binom{j}{2}-kj}\qbinom{k}{j}{q}f_j(q^{-1})
$$
and
$$
\frac{(q^{-r};q^{-d})_k(q^{r-d};q^{-d})_k}{(q^{-d};q^{-d})_k^2}=q^{dk}\cdot
\frac{(q^{r};q^{d})_k(q^{d-r};q^{d})_k}{(q^{d};q^{d})_k^2}.
$$
We have
\begin{align*}
&\sum_{k=0}^{n-1}
\frac{(q^{r};q^{d})_k(q^{d-r};q^{d})_k}{(q^{d};q^{d})_k^2}\cdot g_k(q^{d})\\
\equiv&
(-1)^aq^{d\binom{a+1}2+\frac{(ad+r)(n-1-2a)}{2}}\sum_{k=0}^{n-1}\frac{(q^{r};q^{d})_k(q^{d-r};q^{d})_k}{(q^{d};q^{d})_k^2}\cdot \tilde{g}_k(q^{d})\pmod{\Phi_n(q)^2},
\end{align*}
where
$$
\tilde{g}_k(q)=\sum_{j=0}^k(-1)^jq^{\binom{j}{2}-kj}\qbinom{k}{j}{q}g_j(q).
$$
The similar discussion is also valid when $n$ is even. Thus
\begin{theorem}\label{qsymthm2} Let $n\geq 2$ and $d\geq 1$ with $(n,d)=1$. 
Suppose that $f_1(q),\ldots,f_n(q)\in\Z[q]$ and $r$ is an integer.
Then if $n$ is odd,
\begin{align}\label{qsymodd}
&\sum_{k=0}^{n-1}
\frac{(q^r;q^d)_k(q^{d-r};q^d)_k}{(q^d;q^d)_k^2}\cdot f_k(q^d)\notag\\
\equiv&
(-1)^aq^{d\binom{a+1}2+\frac{(ad+r)(n-1-2a)}{2}}\sum_{k=0}^{n-1}\frac{(q^r;q^d)_k(q^{d-r};q^d)_k}{(q^d;q^d)_k^2}\cdot \tilde{f}_k(q^d)\pmod{\Phi_n(q)^2},
\end{align}
where $a=\langle-r/d\rangle_n$ and
$\Phi_n$ is the $n$-th cyclotomic polynomial.
And when $n$ is even,
\begin{align}\label{qsymeven}
&\sum_{k=0}^{n-1}
\frac{(q^r;q^d)_k(q^{d-r};q^d)_k}{(q^d;q^d)_k^2}\cdot f_k(q^d)\notag\\
\equiv&
(-1)^{a+\frac{ad+r}{n}}q^{d\binom{a+1}2+\frac{(ad+r)(n-1-2a)}{2}}\sum_{k=0}^{n-1}\frac{(q^r;q^d)_k(q^{d-r};q^d)_k}{(q^d;q^d)_k^2}\cdot \tilde{f}_k(q^d)\pmod{\Phi_n(q)^2}.
\end{align}
\end{theorem}

\section{Proof of Theorem \ref{qsymthm}}
\setcounter{theorem}{0}\setcounter{lemma}{0}\setcounter{equation}{0}

Below we slightly extend the notion of congruence. Let
$$
{\mathfrak X}(q)=\Big\{\sum_{j=1}^sc_jq^{\beta_j}:\,s\geq 1,\ c_1,\ldots,c_s\in\C,\ \beta_1,\ldots,\beta_s\in\R\Big\}.
$$
For $f(q),g(q),h(q)\in{\mathfrak X}(q)$, we say $$
f(q)\equiv g(q)\pmod{h(q)},$$ provided
$$
\frac{f(q)-g(q)}{h(q)}\in{\mathfrak X}(q).
$$
\begin{theorem}
Suppose that $n\geq 2$, $\alpha\in\Q$ and the denominator of $\alpha$ is prime to $n$. Then when $n$ is odd,
\begin{align*}
&(-1)^aq^{\binom{a+1}2+sna-s\binom{n}{2}}\sum_{k=0}^{n-1}q^{k^2+k}\qbinom{\alpha}{k}q\qbinom{-1-\alpha}{k}qf_k(q)\\
\equiv&
\sum_{k=0}^{n-1}q^{k^2+k}\qbinom{\alpha}{k}q\qbinom{-1-\alpha}{k}q\hat{f}_k(q)\pmod{\Phi_n(q)^2},
\end{align*}
where $a=\langle\alpha\rangle_n$ and
$s=(\alpha-a)/n$. And if $n$ is even, then
\begin{align*}
&(-1)^{a+s}q^{\binom{a+1}2+sna-s\binom{n}{2}}\sum_{k=0}^{n-1}q^{k^2+k}\qbinom{\alpha}{k}q\qbinom{-1-\alpha}{k}qf_k(q)\\
\equiv&
\sum_{k=0}^{n-1}q^{k^2+k}\qbinom{\alpha}{k}q\qbinom{-1-\alpha}{k}q\hat{f}_k(q)\pmod{\Phi_n(q)^2}.
\end{align*}

\end{theorem}

Let
$$
\hat{f}_k(q)=\sum_{i=0}^k(-1)^iq^{\binom{i+1}2}\qbinom{k}{i}qf_i(q).
$$
Suppose that $p$ is an odd prime and $\alpha$ are positive integer,$0\leq <\alpha>_{p} \leq p-1$.
Let $a=\langle\alpha\rangle_n$.
Then
\begin{align*}
\end{align*}
\begin{proof}
Using the $q$-Chu-Vandemonde identity, we have
\begin{align*}
\qbinom{a+sn}{k}{q} &=\sum_{j=0}^{k} q^{(sn-j)(k-j)}\qbinom{sn}{j}{q}\qbinom{a}{k-j}{q}\\
&\equiv q^{snk} \qbinom{a}{k}{q}+\sum_{j=1}^{k}q^{-j(k-j)}\qbinom{sn}{j}{q}\qbinom{a}{k-j}{q}\pmod{\Phi_n(q)^2}.
\end{align*}
Notice that
$$\qbinom{sn}{j}{q}=\frac{[sn]_q}{[j]_q}\qbinom{sn-1}{j-1}{q}\equiv 0\pmod{\Phi_n(q)}$$
and
$$ \qbinom{sn-1}{j-1}q\equiv\qbinom{-1}{j-1}q=(-1)^{j-1} q^{-\binom{j}{2}}\pmod{\Phi_n(q)}.
$$
So
\begin{align*}
\qbinom{a+sn}{k}q\equiv q^{snk}\qbinom{a}{k}q+\sum_{j=1}^k\frac{[sn]_q}{[j]_q}\cdot(-1)^{j-1}q^{-\binom{j}{2}-j(k-j)}\qbinom{a}{k-j}q\pmod{\Phi_n(q)^2}.
\end{align*}
Similarly, we also have
\begin{align*}
&\qbinom{-1-a-sn}{k}q\\
\equiv&q^{-snk}\qbinom{-1-a}{k}q+\sum_{j=1}^k\frac{[-sn]_q}{[j]_q}\cdot(-1)^{j-1}q^{-\binom{j}{2}-j(k-j)}\qbinom{-1-a}{k-j}q\pmod{\Phi_n(q)^2}.
\end{align*}
Thus
\begin{align}
&\sum_{k=0}^{n-1}q^{k^2+k}\qbinom{\alpha}{k}q\qbinom{-1-\alpha}{k}q((-1)^aq^{\binom{a+1}2}f_k(q)-\hat{f}_k(q))\notag\\
\equiv
&\sum_{k=0}^{n-1}q^{k^2+k}\qbinom{a}{k}q\qbinom{-1-a}{k}q((-1)^aq^{\binom{a+1}2}f_k(q)-\hat{f}_k(q))\label{part1}\\
&-s\sum_{k=1}^{n-1}q^{k^2+k}\qbinom{a}{k}qf_k(q)\sum_{j=1}^k(-1)^{a+j-1}q^{\binom{a+1}2-\binom{j}{2}-j(k-j)}\cdot\frac{[n]_q}{[j]_q}\qbinom{-1-a}{k-j}q\label{part2}\\
&+s\sum_{k=1}^{n-1}q^{k^2+k}\qbinom{a}{k}q\hat{f}_k(q)\sum_{j=1}^k(-1)^{j-1}q^{-\binom{j}{2}-j(k-j)}\cdot\frac{[n]_q}{[j]_q}\qbinom{-1-a}{k-j}q\label{part3}\\
&+s\sum_{k=1}^{n-1}q^{k^2+k}\qbinom{-1-a}{k}qf_k(q)\sum_{j=1}^k(-1)^{a+j-1}q^{\binom{a+1}2-\binom{j}{2}-j(k-j)}\cdot\frac{[n]_q}{[j]_q}\qbinom{a}{k-j}q\label{part4}\\
&-s\sum_{k=1}^{n-1}q^{k^2+k}\qbinom{-1-a}{k}q\hat{f}_k(q)\sum_{j=1}^k(-1)^{j-1}q^{-\binom{j}{2}-j(k-j)}\cdot\frac{[n]_q}{[j]_q}\qbinom{a}{k-j}q\pmod{\Phi_n(q)^2}.\label{part5}
\end{align}
First, we consider (\ref{part1}). Clearly
\begin{align*}
&\sum_{k=0}^{n-1} q^{k^2+k}\qbinom{a}{k}{q}\qbinom{-1-a}{k}{q}\hat{f}_k(q)\\
=&\sum_{k=0}^{a} q^{k^2+k}\qbinom{a}{k}{q}\qbinom{-1-a}{k}{q} \sum_{j=0}^{k}\qbinom{k}{j}{q}(-1)^{j}q^{\binom{j+1}{2}}f_j(q)\\
=&\sum_{j=0}^{a}(-1)^j q^{\binom{j+1}2}f_j(q)\sum_{k=j}^{a}q^{k^2+k}\qbinom{a}{k}q\qbinom{-1-a}{j}q\qbinom{-1-a-j}{k-j}q\\
=&\sum_{j=0}^{a}(-1)^j q^{\binom{j+1}2}f_j(q)\qbinom{-1-a}{j}{q}q^{a^2+a}\sum_{k=j}^{a}q^{(a-k)(-1-a-k)}\qbinom{a}{a-k}q\qbinom{-1-a-j}{k-j}q\\
=&\sum_{j=0}^{a}(-1)^j q^{\binom{j+1}2}f_j(q)\qbinom{-1-a}{j}{q} q^{a^2+a}\qbinom{-1-j}{a-j}q\\
=&\sum_{j=0}^{a}(-1)^a q^{j^2+j+\binom{a+1}{2}}f_j(q)\qbinom{-1-a}{j}{q} \qbinom{a}{j}{q},
\end{align*}
where we use the $q$-Chu-Vandemonde identity in the fourth equality.
So (\ref{part1}) always vanishes.

Note that
\begin{align*}
&\sum_{j=1}^{k}(-1)^{j-1}q^{-\binom{j}{2}-j(k-j)}\cdot\frac{[n]_q}{[j]_q}\qbinom{-1-a}{k-j}{q}
\equiv\sum_{j=1}^{k}\qbinom{n-1}{j-1}{q}\cdot q^{-j(k-j)}\cdot\frac{[n]_q}{[j]_q}\qbinom{-1-a}{k-j}{q}\\
\equiv&\sum_{j=1}^{k}q^{(n-j)(k-j)}\qbinom{n}{j}{q}\qbinom{-1-a}{k-j}{q}
=\qbinom{n-1-a}{k}{q}-q^{nk}\qbinom{-1-a}{k}{q}\pmod{\Phi_n(q)^2}.
\end{align*}
Hence
\begin{align}
(\ref{part2})\equiv&-s(-1)^{a}q^{\binom{a+1}{2}}\sum_{k=1}^{n-1}q^{k^2+k}\qbinom{a}{k}{q}f_k(q)\cdot\bigg(\qbinom{n-1-a}{k}{q}-q^{nk}\qbinom{-1-a}{k}{q}\bigg)\\
\equiv&-s(-1)^{a}q^{\binom{a+1}{2}}\sum_{k=0}^{n-1}q^{k^2+k}\qbinom{a}{k}{q}\qbinom{n-1-a}{k}{q}f_k(q)\label{part21}\\
&+s(-1)^{a}q^{\binom{a+1}{2}}\sum_{k=0}^{n-1} q^{k^2+k+nk}\qbinom{a}{k}{q}\qbinom{-1-a}{k}{q}f_k(q)\pmod{\Phi_n(q)^2}.\label{part22}
\end{align}
Similarly, we have
\begin{align*}
&\sum_{j=1}^{k}(-1)^{j-1}q^{-\binom{j}{2}-j(k-j)}\frac{[n]_q}{[j]_q}\qbinom{a-n}{k-j}{q}\\
\equiv&\sum_{j=1}^{k}q^{(n-j)(k-j)}\qbinom{n}{j}q\qbinom{a-n}{k-j}{q}=\qbinom{a}{k}{q}-q^{nk}\qbinom{a-n}{k}{q}\pmod{\Phi_n(q)^2}.\end{align*}
So
\begin{align}
(\ref{part4})\equiv&s\sum_{k=1}^{n-1-a} q^{k^2+k+\binom{a+1}{2}}\qbinom{n-1-a}{k}{q}f_k(q)\sum_{j=1}^{k}(-1)^{a+j-1}q^{-\binom{j}{2}-j(k-j)}\cdot\frac{[n]_q}{[j]_q}\qbinom{a-n}{k-j}{q}\notag\\
\equiv&s(-1)^{a}q^{\binom{a+1}{2}}\sum_{k=0}^{n-1-a} q^{k^2+k}\qbinom{n-1-a}{k}{q}\qbinom{a}{k}{q}f_k(q)\label{part41}\\
&-s(-1)^{a}q^{\binom{a+1}{2}}\sum_{k=0}^{n-1-a} q^{k^2+k+nk}\qbinom{n-1-a}{k}{q}\qbinom{a-n}{k}{q}f_k(q)\pmod{\Phi_p(q)^2}.\label{part42}
\end{align}

On the other hand, clearly
\begin{align*}
&\sum_{k=1}^{n-1}q^{k^2+k}\qbinom{a}{k}q\hat{f}_k(q)\sum_{j=1}^k(-1)^{j-1}q^{-\binom{j}{2}-j(k-j)}\cdot\frac{[n]_q}{[j]_q}\qbinom{-1-a}{k-j}q\\
=&\sum_{k=1}^{n-1}q^{k^2+k}\qbinom{a}{k}q\sum_{i=0}^k(-1)^iq^{\binom{i+1}{2}}\qbinom{k}{i}qf_i(q)\sum_{j=1}^k(-1)^{j-1}q^{-\binom{j}{2}-j(k-j)}\cdot\frac{[n]_q}{[j]_q}\qbinom{-1-a}{k-j}q\\
=&\sum_{i=0}^a(-1)^iq^{\binom{i+1}{2}}\qbinom{a}{i}qf_i(q)\sum_{j=1}^{a}\frac{(-1)^{j-1}q^{-\binom{j}{2}}[n]_q}{[j]_q}\sum_{k=j}^{a}q^{k^2+k-j(k-j)}\qbinom{a-i}{a-k}q\qbinom{-1-a}{k-j}q\\
=&\sum_{i=0}^a(-1)^iq^{\binom{i+1}{2}}\qbinom{a}{i}qf_i(q)\sum_{j=1}^{a}\frac{(-1)^{j-1}q^{-\binom{j}{2}}[n]_q}{[j]_q}
\cdot q^{a^2+a+j^2-aj}\qbinom{-1-i}{a-j}q,
\end{align*}
where in the last step we use the $q$-Chu-Vandemonde identity. Thus
\begin{align}
(\ref{part3})=&s\sum_{i=0}^a(-1)^iq^{\binom{i+1}{2}}\qbinom{a}{i}qf_i(q)\sum_{j=1}^{a}\frac{(-1)^{j-1}q^{a^2+a+j^2-aj-\binom{j}{2}}[n]_q}{[j]_q}\qbinom{-1-i}{a-j}q\notag\\
\equiv&s\sum_{i=0}^a(-1)^iq^{\binom{i+1}{2}}\qbinom{a}{i}qf_i(q)\sum_{j=1}^{a}q^{a^2+a+j^2-aj}\qbinom{n}{j}q\qbinom{-1-i}{a-j}q\notag\\
\equiv&sq^{a^2+a}\sum_{i=0}^a(-1)^iq^{\binom{i+1}{2}}\qbinom{a}{i}qf_i(q)\sum_{j=1}^{a}q^{(n-j)(a-j)}\qbinom{n}{j}q\qbinom{-1-i}{a-j}q\notag\\
=&sq^{a^2+a}\sum_{i=0}^a(-1)^iq^{\binom{i+1}{2}}\qbinom{a}{i}qf_i(q)
\cdot\qbinom{n-1-i}{a}q\label{part31}\\
-&sq^{a^2+a}\sum_{i=0}^a(-1)^iq^{\binom{i+1}{2}}\qbinom{a}{i}qf_i(q)
\cdot q^{na}\qbinom{-1-i}{a}q\pmod{\Phi_n(q)^2}.\label{part32}
\end{align}
Similarly, noting that
\begin{align*}
&\sum_{k=j}^{n-1-a}q^{k^2+k-j(k-j)}\qbinom{n-1-a-i}{k-i}q\qbinom{a}{k-j}q\\
\equiv&\sum_{k=j}^{n-1-a}q^{a^2+a+j^2+j+aj+(n-1-a-k)(a-k+j)}\qbinom{n-1-a-i}{n-1-a-k}q\qbinom{a}{k-j}q\\
=&q^{a^2+a+j^2+j+aj}\qbinom{n-1-i}{n-1-a-j}q\pmod{\Phi_n(q)},
\end{align*}
we can get
\begin{align}
(\ref{part5})\equiv&-s\sum_{k=1}^{n-1-a}\qbinom{n-1-a}{k}q\hat{f}_k(q)\sum_{j=1}^k\frac{(-1)^{j-1}q^{-\binom{j}{2}}[n]_q}{[j]_q}\cdot q^{k^2+k-j(k-j)}\qbinom{a}{k-j}q\notag\\
\equiv&-s\sum_{i=0}^{n-1-a}(-1)^iq^{\binom{i+1}{2}}\qbinom{n-1-a}{i}qf_i(q)\sum_{j=1}^{n-1-a}\qbinom{n}{j}q\cdot q^{a^2+a+j^2+j+aj}\qbinom{n-1-i}{n-1-a-j}q\notag\\
\equiv&-s\sum_{i=0}^{n-1-a}(-1)^iq^{\binom{i+1}{2}}\qbinom{n-1-a}{i}qf_i(q)\sum_{j=1}^{n-1-a}q^{a^2+a+(n-j)(n-1-a-j)}\qbinom{n}{j}q\qbinom{-1-i}{n-1-a-j}q\notag\\
=&-s\sum_{i=0}^{n-1-a}(-1)^iq^{a^2+a+\binom{i+1}{2}}\qbinom{n-1-a}{i}qf_i(q)\cdot\qbinom{n-1-i}{n-1-a}q\label{part51}\\
+&s\sum_{i=0}^{n-1-a}(-1)^iq^{a^2+a+\binom{i+1}{2}}\qbinom{n-1-a}{i}qf_i(q)\cdot q^{n(n-1-a)}\qbinom{-1-i}{n-1-a}q\pmod{\Phi_n(q)^2}.\label{part52}
\end{align}

Clearly
\begin{equation}\label{final1}
(\ref{part21})+(\ref{part41})=0.
\end{equation}
And we also get
\begin{equation}\label{final2}
(\ref{part31})+(\ref{part51})=0,
\end{equation}
since
\begin{align*}
\qbinom{a}{i}q\qbinom{n-1-i}{a}q=&
\qbinom{n-1-i}{i}q\qbinom{n-1-2i}{a-i}q\\
=&
\qbinom{n-1-i}{i}q\qbinom{n-1-2i}{n-1-a-i}q=
\qbinom{n-1-a}{i}q\qbinom{n-1-i}{n-1-a}q.
\end{align*}
Noting that
$$
(-1)^{a}q^{ak+\binom{a+1}{2}}\qbinom{-1-k}{a}q=\qbinom{a+k}{a}q
=\qbinom{a+k}{k}q=
(-1)^{k}q^{ak+\binom{k+1}{2}}\qbinom{-1-a}{k}q,
$$
we have
\begin{equation}\label{final3}
(\ref{part22})+(\ref{part32})=
s(-1)^{a}q^{\binom{a+1}{2}}\sum_{k=0}^{n-1} q^{k^2+k}\qbinom{a}{k}{q}\qbinom{-1-a}{k}{q}f_k(q)\cdot(q^{nk}-q^{na}).
\end{equation}
Similarly, from
$$
(-1)^{n-1-a}q^{\binom{n-a}{2}}\qbinom{-1-k}{n-1-a}q
=(-1)^{k} q^{\binom{k+1}{2}}\qbinom{a-n}{k}q,
$$
it follows that
\begin{align}
&(\ref{part42})+(\ref{part52})\notag\\
=&
s(-1)^{a}q^{\binom{a+1}{2}}\sum_{k=0}^{n-1-a} q^{k^2+k}\qbinom{n-1-a}{k}{q}\qbinom{a-n}{k}{q}f_k(q)
\cdot((-1)^{n-1}q^{\binom{n}{2}}-q^{nk}).\label{final4}
\end{align}
Noting that
$$
1+q^{\frac n2}=\frac{1-q^n}{1-q^{\frac n2}}\equiv 0\pmod{\Phi_n(q)}
$$
for those even $n$,
we always have
$$
(-1)^{n-1}q^{\binom{n}{2}}\equiv1\pmod{\Phi_n(q)}.
$$
Thus by (\ref{final1}), (\ref{final2}), (\ref{final3}) and (\ref{final4}), we get
\begin{align*}
&\sum_{k=0}^{n-1}q^{k^2+k}\qbinom{\alpha}{k}q\qbinom{-1-\alpha}{k}q((-1)^aq^{\binom{a+1}2}f_k(q)-\hat{f}_k(q))
\\\equiv&
(-1)^{a}sq^{\binom{a+1}{2}}\sum_{k=0}^{n-1}q^{k^2+k}((-1)^{n-1}q^{\binom{n}{2}}-q^{na})\qbinom{-1-a}{k}{q}\qbinom{a}{k}{q}f_k(q)\\
\equiv&
(-1)^{a}sq^{\binom{a+1}{2}}\sum_{k=0}^{n-1}q^{k^2+k}(1-(-1)^{n-1}q^{na-\binom{n}{2}})\qbinom{-1-\alpha}{k}{q}\qbinom{\alpha}{k}{q}f_k(q)\pmod{\Phi_n(q)^2}.
\end{align*}

When $n$ is odd, we have
$$
s(1-q^{na-s\binom{n}{2}})\equiv
1-q^{sna-s\binom{n}{2}}\pmod{\Phi_n(q)^2}.
$$
So
\begin{align*}
&(-1)^aq^{\binom{a+1}2+sna-s\binom{n}{2}}\sum_{k=0}^{n-1}q^{k^2+k}\qbinom{\alpha}{k}q\qbinom{-1-\alpha}{k}qf_k(q)\\
\equiv&
\sum_{k=0}^{n-1}q^{k^2+k}\qbinom{\alpha}{k}q\qbinom{-1-\alpha}{k}q\hat{f}_k(q)\pmod{\Phi_n(q)^2}.
\end{align*}
Suppose that $n$ is even. Then
$$
s(1-(-1)^{n-1}q^{na-s\binom{n}{2}})\equiv
1-(-1)^sq^{sna-s\binom{n}{2}}\pmod{\Phi_n(q)^2}.
$$
Then we have
\begin{align*}
&(-1)^{a+s}q^{\binom{a+1}2+sna-s\binom{n}{2}}\sum_{k=0}^{n-1}q^{k^2+k}\qbinom{\alpha}{k}q\qbinom{-1-\alpha}{k}qf_k(q)\\
\equiv&
\sum_{k=0}^{n-1}q^{k^2+k}\qbinom{\alpha}{k}q\qbinom{-1-\alpha}{k}q\hat{f}_k(q)\pmod{\Phi_n(q)^2}.
\end{align*}
\end{proof}

\begin{Rem}
Let 
$$
f_k(q)=x^k.
$$
Then
$$
\hat{f}_k(q)=\sum_{j=0}^k(-1)^jq^{\binom{j+1}{2}}\qbinom{k}{j}qx^j=(xq;q)_k.
$$
Then Theorem \ref{qsymthm} implies a conjecture of Guo and Zeng \cite[Conjecture 7.1]{GZ14}.
\end{Rem}
\begin{Rem}
Let 
$$
f_k(q)=\frac{q^k(x;q)_k}{(q;q)_k}.
$$
Then it is not difficult to prove that
$$
\tilde{f}_k(q)=\sum_{j=0}^k(-1)^jq^{\binom{j}{2}-kj}\qbinom{k}{j}q\cdot\frac{q^j(x;q)_j}{(q;q)_j}=
\frac{(xq^{-1};q^{-1})_k}{(q^{-1};q^{-1})_k}.
$$
Define
$$
\cP_n(\alpha,x;q)=\sum_{k=0}^{n-1}q^{k^2+k}\qbinom{\alpha}{k}q\qbinom{-1-\alpha}{k}q\cdot
\frac{(x;q)_k}{(q;q)_k}.
$$
Then in view of Theorem \ref{qsymthm2}, we can get
$$
\cP_n(-r/d,x;q^d)\equiv
(-1)^aq^{d\binom{a+1}2+\frac{(ad+r)(n-1-2a)}{2}}\cP_n(-r/d,xq^{-d};q^{-d})\pmod{\Phi_n(q)^2}
$$
for odd $n\geq 3$,
where $a=\langle-r/d\rangle_n$. This is a $q$-analogue of \cite[Theorem 2.2]{Sun14}.
\end{Rem}

\end{document}